\documentclass[11pt]{article}
\usepackage{osa2}
\usepackage{overcite}  

\begin{document}


\title{Phase retrieval by iterated projections}

\author{Veit Elser}
\address{Department of Physics, Cornell University, Ithaca NY 14853-2501}
\email{ve10@cornell.edu}

\begin{abstract}
Several strategies in phase retrieval are unified by an iterative ``difference map" constructed from a pair of elementary projections and a single real parameter $\beta$. For the standard application in optics, where the two projections implement Fourier modulus and object support constraints respectively, the difference map reproduces the ``hybrid" form of Fienup's \cite{Fienup1} input-output map for $\beta = 1$. Other values of $\beta$ are equally effective in retrieving phases but have no input-output counterparts. The geometric construction of the difference map illuminates the distinction between its fixed points and the recovered object, as well as the mechanism whereby stagnation is avoided. When support constraints are replaced by object histogram or atomicity constraints, the difference map lends itself to crystallographic phase retrieval. Numerical experiments with synthetic data suggest that structures with hundreds of atoms can be solved.
\end{abstract}

\section{Introduction}
In the standard phase retrieval problem one has available the modulus of the Fourier transform of an ``object", usually in two or three dimensions, and a set of {\it a priori} constraints that a proper Fourier reconstruction of the object must satisfy. Usually one assumes that, up to translation and inversion of the object, the reconstruction is unique. In other words, there is essentially a unique set of phases that, when combined with the known Fourier modulus data, satisfy all the {\it a priori} constraints. It is in this sense that the unknown phases are said to be ``retrieved" from the modulus data.

Currently there are no practical algorithms for phase retrieval, that is, procedures where a solution is guaranteed at a computational cost that grows modestly with the size of the problem. Rather, a number of heuristic strategies have been developed that in many situations perform very well. A dominant theme of these practical methods is an iteration scheme which modifies the object in each cycle until both the Fourier moduli are correct and all the {\it a priori} constraints are satisfied. Although nothing has been proven about the rate of convergence of such schemes, the body of favorable empirical evidence in both imaging\cite{Stark1} and crystallographic\cite{directMethods} applications is substantial.

Expressed in the most general terms, phase retrieval is the problem of finding an element in a large set that simultaneously has properties $A$, $B$, etc., given only the ability to find elements having these properties separately. If property $A$ corresponds to the object having the correct Fourier modulus, property $B$ might represent, for example, a support constraint. While it is very easy to construct objects having the correct Fourier modulus, or a prescribed support, there are no known methods for directly constructing an object (if one exists) that has both properties. Accepting this basic limitation, we can nevertheless attempt to build the desired construction using the simple constructions for properties $A$, $B$, etc. as {\it elementary operations}.

A natural setting\cite{Stark2} for reconstructing an object is the $N$ dimensional Euclidean vector space $E^N$, whose components correspond to the values of $N$ pixels, in an imaging application, or in the case of crystallography, a sampling of the electron density at $N$ regular grid points in the unit cell. The main benefit of the Euclidean space is the ability to define a distance between objects. It is convenient to consider a complex Euclidean space, where unitary transformations such as the Fourier transform may act. Since distance is preserved by unitary transformations, distances can be measured in either the object or Fourier domains.

With the ability to measure distance we can define more precisely the elementary operations for constructing objects with properties $A$, $B$, etc.\cite{Stark2} These are the ``projections", which begin with some point $\rho \in E^N$, and find a point $\pi_{A}(\rho)$ having property $A$, for example, and whose distance from $\rho$ is a minimum. For most of the projections usually considered in imaging and crystallography, the point with minimum distance is essentially unique and can be computed in a time that grows no faster than $N \log N$. This is important, since we anticipate that many elementary projections will be required to arrive at the desired object having all the properties $A$, $B$, etc.

The earliest application of projections to reconstruct an image from Fourier modulus and support constraints is the Gerchberg-Saxton\cite{GerchSaxt} map: $G = \pi_{\rm supp} \circ \pi_{\rm mod}$. When applied to an initial randomly chosen point, the idea is that successive iterates might simultaneously get closer, both to the subspace of $E^N$ having the correct Fourier modulus, as well as the subspace corresponding to the support constraint. Owing to the non-convexity\cite{Stark3} of the subspace of correct Fourier modulus, this map has the problem that it can encounter fixed points which have the correct support but fail to have the correct modulus. The Gerchberg-Saxton map in this situation moves the point from one constraint subspace to the other, and then back to the original point; a point common to both subspaces has not been identified. This phenomenon, potentially present for all projections where the corresponding subspace is non-convex, is known as ``stagnation".

To avoid stagnation, various generalizations\cite{Fienup1, Stark3} of the Gerchberg-Saxton map have been proposed. In imaging, the most successful of these have been schemes where the iterates are not confined to the constraint subspaces. In particular, objects realized as general linear combinations of other objects are generated that take advantage of the vector space nature of the search space. For these generalized maps it is usually the case that none of the iterates enroute to the solution satisfy {\it any} of the constraints. Fienup\cite{Fienup1} has devised a family of maps motivated by ideas from control theory, of which the ``hybrid input-output" map stands out as the most successful in a wide variety of imaging applications. The theoretical understanding of the hybrid input-output map, however, is very incomplete. There is not even a compelling argument that sets the ``hybrid" version of the map apart from its less successful input-output relatives. Also, whereas the hybrid input-output map was developed for imaging applications with a support constraint, one would like to be able to exploit its excellent performance characteristics in applications such as crystallography, where support constraints are usually not available.

This paper describes a general map formed by taking the difference of a pair of elementary projections. The advantage of such a ``difference map" is that stagnation, in a strict sense, can be ruled out. By optimizing local convergence properties, the detailed form of the difference map is determined up to a single real parameter $\beta$. Interestingly, for the special case of support constraints, the value $\beta = 1$ reproduces an input-output type of map, and it is exactly of the hybrid variety. Other values of $\beta$, such as $-1$, give maps for which there are no input-output counterparts and yet whose performance in numerical experiments is also very good. Perhaps most significant is the versatility of the difference map with respect to the elementary projections from which it is built. In particular, all of the examples reconstructed in numerical experiments with support projection can, when implemented by the difference map, be reconstructed using an elementary projection onto a known distribution of object values. The set of object values, or histogram, is a form of {\it a priori} constraint that can be exploited in crystallography. The difference map thus provides a link between the most successful phase retrieval method in imaging, the hybrid input-output map, and the crystallographic phase problem. 

\section{Examples of elementary projections}

A projection, of arbitrary points in the $N$ dimensional Euclidean space of objects, to points on a subspace, representing a particular constraint, is said to be {\it elementary} if it is (i) distance minimizing and unique, and (ii) easy to compute. In imaging and crystallography, ``easy to compute" translates to the statement that the operation count scale essentially as $N$. Below is a collection of elementary projections that have been employed in various phase retrieval schemes:

\noindent {\bf Fourier modulus.} It is convenient to express all projections with respect to a common basis in $E^N$, which we take to be the object domain. Because the Fourier modulus projection is more naturally expressed in the Fourier domain, we write
\begin{equation}
\pi_{\rm mod} = \mathcal{F}^{-1}\cdot \tilde{\pi}_{\rm mod}\cdot \mathcal{F}\quad ,
\end{equation}
where $\mathcal{F}$ is the unitary transformation to the Fourier domain, and $\tilde{\pi}_{\rm mod}$ is the projection operator which acts componentwise, that is, on each pixel of the (discrete) Fourier transform. Geometrically, each complex component in the Fourier transform is mapped by $\tilde{\pi}_{\rm mod}$ to the nearest point on a circle having a prescribed radius (the Fourier modulus). Thus it follows that $\pi_{\rm mod}$ is distance minimizing and essentially unique (the exception being the measure-zero set of complex numbers with zero modulus). Since the computation of $\pi_{\rm mod}$ is linear in $N$, and the multiplication by $\mathcal{F}$ and $\mathcal{F}^{-1}$ requires $N \log N$ operations when using the FFT algorithm, Fourier modulus projection is easy to compute. Since circles are not convex, the subspace of objects satisfying the Fourier modulus constraints, $C_{\rm mod}$, is also non-convex.

\noindent {\bf Support.} There is a unique, distance minimizing map onto the subspace $C_{\rm supp}$ of objects having a specified support $S$. If $\rho_{n}$ is the value of pixel $n$ in the object domain, then support projection is the map
\begin{equation}
\pi_{\rm supp}\; \colon\; \rho_{n} \mapsto \rho_{n}^{\prime} = \left\{
\begin{array}{ll}
\rho_{n} & \mbox{if $n \in S$} \\
0 & \mbox{if $n \notin S$.} 
\end{array}
\right.
\end{equation}
Since $C_{\rm supp}$ is a linear subspace it is convex; projection computations require at most $N$ operations.

\noindent {\bf Positivity.} For real-valued objects one can impose positivity. The unique, distance minimizing map to the corresponding subspace $C_{\rm pos}$, sets all negative pixels values to zero and leaves the positive pixel values unchanged. Like $C_{\rm supp}$, $C_{\rm pos}$ is convex. Moreover, the projections $\pi_{\rm supp}$ and $\pi_{\rm pos}$ commute and can be combined into a single elementary projection.

\noindent {\bf Histogram.} Histogram projection\cite{histMatch} is complementary to support projection, if we regard the object as a function from the support $S$ to a set of values, or ``histogram", $H$. Here we define the object histogram $H$ to be the set of $N$ pixel values without regard to pixel position. For real-valued objects a distance minimizing map $\pi_{\rm hist}$, onto the subspace of objects having histogram $H$, is unique and easily computed. One begins by sorting $H = \{h_1, h_2, \dots \}$ and also finding an ordering of the pixels, $n \mapsto o(n)$, such that their values, $\{\rho_{o(1)}, \rho_{o(2)}, \dots \}$ are also sorted. Histogram projection is then given by the map
\begin{equation}
\pi_{\rm hist}\; \colon \; \rho_{o(n)} \mapsto h_{n} .
\end{equation}
It is straightforward to show that $\pi_{\rm hist}$ is distance minimizing. Since $N \log N$ operations are required to sort $N$ real numbers, histogram projection is easy to compute. For complex valued objects it is probably not possible to compute the projection to a specified complex histogram in order $N$ operations, although approximate distance minimizing maps can be computed with this effort\cite{multiDimHist}. The subspace of objects having a prescribed histogram, $C_{\rm hist}$, is the point set formed by applying all permutations to the point $\{h_1, h_2, \dots \}$, and as such is non-convex.

\noindent {\bf Atomicity.} The projection of greatest relevance to crystallography (and possibly astronomy) is the distance minimizing map to the set of objects consisting of a known number $M$ of non-overlapping atoms (or stars). It is not necessary to make the restriction to equal atoms, although this simplifies the computation of the projection. Our model of atoms will be objects of small support, usually just a $3 \times 3$ array of pixels ($3\times 3\times 3$ in three dimensions). The object values for each type of atom are also specified, and allowance is made for the possibility that the actual center of the atom may be located between the pixel centers (see Appendix for details). A unique, distance minimizing map is most easily constructed for the case of identical atoms with support on a single pixel. This situation may equivalently be characterized by its histogram: a set of $M$ identical positive values $\rho_{+}$, representing the atoms, and $N-M$ zeros. Atom projection, in this case, corresponds simply to histogram projection: the $M$ largest pixels in the object are set to $\rho_{+}$ and the remainder set to zero. The difference between histogram and atom projection emerges when the support of each atom is larger than a single pixel; a graphical illustration is given in Figure 1.
\begin{figure}[h]
\centerline{\scalebox{1.}{\includegraphics{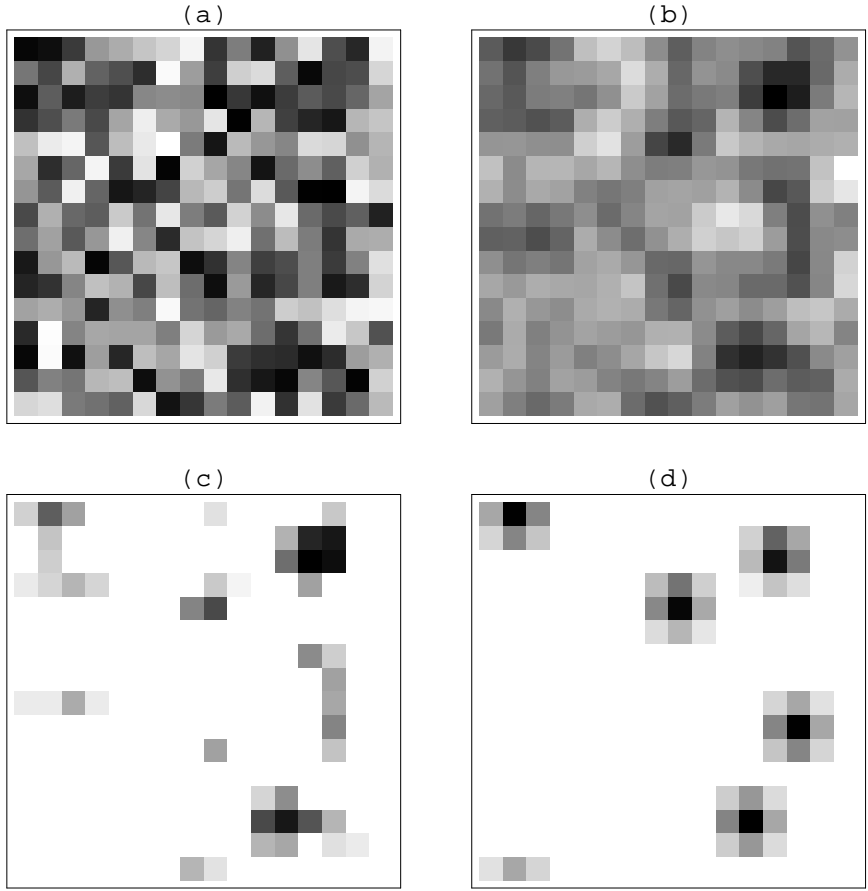}}}
\caption{Examples of elementary projections. (a) Random positive density; (b) Fourier modulus projection of (a); (c) histogram projection of (b); (d) atomicity projection of (b). 
}
\end{figure}
In general terms we can say that, whereas histogram projection completely neglects spatial correlations of the pixel values in the object domain, with atom projection one at least incorporates short range correlations.

It is more difficult to construct atom projections for multi-pixel atoms, and consequently, to rigorously prove minimality. In practice it probably suffices to map the object reasonably near the true distance minimizing point on the subspace of atomic objects, $C_{\rm atom}$. The algorithm used in the numerical experiments described in Section 7 begins by sorting the pixels by value and, beginning with the largest, associates them with the specified number of atoms while observing the restriction that each new atom's support does not overlap the support of previously identified atoms. The mapping is then fine-tuned by computing the centroid of the object values within each atom's support and setting the values on the support equal to that of an atom having that particular fractional location with respect to the pixels. More details are given in the Appendix.

Apart from the minor complication due to non-overlapping supports, the subspace $C_{\rm atom}$ is just the product of tori corresponding to the locations of atoms within the periodic crystal unit cell (a non-convex set). This assumes that atom locations can move continuously between pixels. In practice this is only approximately true, since the support of each atom must move one pixel at a time. The discreteness of pixels also manifests itself in the relationship of $C_{\rm atom}$ to $C_{\rm hist}$, which should be one of containment: $C_{\rm atom} \subset C_{\rm hist}$. For this to be true we see that $C_{\rm hist}$ must contain object values for all fractional displacements of an atom relative to the pixels. Clearly this is not true if we interpret $C_{\rm hist}$ as the point set formed by taking all permutations of a given set of $N$ histogram values. The true histogram should include a continuum of values, corresponding to continuous translations of the atoms. For the practical implementation of histogram projection, however, it suffices to sample just $N$ representative values from the continuous distribution and use these as described above.

Some observations of a general nature can be made if the constraint subspaces are smooth submanifolds of $E^N$. While this is the case for the Fourier modulus and support constraints, it applies approximately for atomicity since $C_{\rm atom}$ has a nearly smooth parametrization by continuous atomic positions on the $d$-dimensional torus. If one is seeking an object $\rho \in E^N$ in the intersection of two smooth constraint spaces, say $C_A$ and $C_B$, then just the dimensionalities of these spaces can provide useful information. For generic embeddings in $E^N$, for example, one would not expect a unique solution if $\dim{C_A}+\dim{C_B} > N$. For phase retrieval ($C_A = C_{\rm mod}$) to be well posed, this implies $\dim{C_B} < N/2$, that is, an object support that is smaller than half of the image, if $C_B = C_{\rm supp}$, or a number of atoms $M$ smaller than $N/(2 d)$, if $C_B = C_{\rm atom}$. The converse of this is that when a generic embedding predicts an empty intersection, the fact that a solution is {\it a priori} known to exist implies it is unique. However, as the case of phase retrieval with support constraints in one dimension\cite{nonunique1d} shows, this form of reckoning can fail. Nevertheless, experiments (Sec. 7) with atomicity and histogram constraints suggest that nonuniqueness in the case of overdetermined constraints is the exception rather than the rule. 

\section{The difference map}

Given two arbitrary projections $\pi_1$ and $\pi_2$, we consider the ``difference map" $D \colon E^N \to E^N$ defined by
\begin{equation}
D = 1 + \beta \Delta \quad,
\end{equation}
where $\beta$ is a nonzero real parameter and
\begin{equation}\label{diffMap}
\Delta = \pi_1 \circ f_2 - \pi_2 \circ f_1
\end{equation}
is the difference of the two projection operators, each composed with a map $f_i \colon E^N \to E^N$. The detailed form of the maps $f_i$ is secondary to the global behavior of the difference map and is discussed in the next Section.  A fixed point of $D$, $\rho^*$, is characterized by $\Delta( \rho^* ) = 0$, or
\begin{equation}\label{projFixedPt}
(\pi_1 \circ f_2) ( \rho^* ) = \rho_{1\cap 2} = (\pi_2 \circ f_1) ( \rho^* ) \quad,
\end{equation}
where $\rho_{1\cap 2}$ is required to lie in the intersection of the corresponding constraint subspaces. If $\pi_2 = \pi_{\rm mod}$, say, and $\pi_1$ represents some object domain constraint, then $\rho_{1\cap 2}$ is a solution of that particular instance of the phase problem. We note that in general $\rho_{1\cap 2} \neq \rho^*$. The standard approach for finding solutions is to begin with a randomly chosen point $\rho(0) \in E^N$, iterate the map $D$ until a fixed point $\rho^*$ is found, and then use (\ref{projFixedPt}) to find $\rho_{1\cap 2}$. Even if the solution $\rho_{1\cap 2}$ is unique (up to translation and inversion), the set of fixed points is in general the large space given by
\begin{equation}\label{fixedPts}
(\pi_1 \circ f_2)^{-1} ( \rho_{1\cap 2} ) \,\cap\, (\pi_2 \circ f_1)^{-1} ( \rho_{1\cap 2} )\quad.
\end{equation}
One practical consequence of this is that the fixed point object $\rho^*$ will appear to be contaminated with noise whose origin, ultimately, is the randomness of the starting point $\rho(0)$. This is illustrated by the numerical experiments in Section 7. 

The progress of the iterates $\rho(i)$ can be monitored by keeping a record of the norm of the differences
\begin{equation}
e_i = \| \Delta( \rho(i) ) \| \quad,
\end{equation}
where $\|\cdot \|$ denotes the Euclidean norm. The ``error" $e_i$ has the geometrical interpretation as the currently achieved distance between the two constraint subspaces, $C_1$ and $C_2$. When this distance becomes sufficiently small, an object has been found that satisfies both sets of constraints. The behavior of $e_i$ is, in general, nonmonotonic.

It is convenient to scale the constraint subspaces so that they lie on the sphere of unit norm. The magnitude of the error estimate is then related to the angular separation of the constraint subspaces by
\begin{equation}
e = 2 \sin\left( {\theta_{1\cap2}}\over{2}\right)\quad.
\end{equation}
While small values of $\theta_{1\cap2}$ imply a solution is nearby, small values encountered during the early iterates are a cause for concern. This is because it is highly improbable that the random starting point $\rho(0)$ is within a small angle of the solution, or more precisely, the set of fixed points (\ref{fixedPts}). Small initial errors (angles) are more likely to be an indication that the object constraints are too weak for reliable phase retrieval, as for example an object support that occupies nearly half the image. 

\section{Local convergence and traps}

The maps $f_i$ in the definition of the difference map (\ref{diffMap}) must be chosen with care if the set of fixed points is to be attractive. For example, choosing the identity map for both $f_1$ and $f_2$ does not give attractive fixed points. In order to limit the possibilities for these maps we use a geometric approach that relies on the elementary projections $\pi_1$ and $\pi_2$. Let $\rho \in E^N$ be the current iterate (object). Using projections we can obtain two additional points, $\pi_1(\rho)$ and $\pi_2(\rho)$. Pairs of points determine lines; in particular,
\begin{equation}\label{fMap}
f_i(\rho) = (1 + \gamma_i) \pi_i(\rho) - \gamma_i \rho
\end{equation}
is a general point, parametrized by a real number $\gamma_i$, on the line defined by $\rho$ and $\pi_i(\rho)$. Using (\ref{fMap}) as a definition of the maps $f_i$ might be viewed as taking the first step beyond simply using identity maps for the $f_i$. Indeed, the parameter choice $\gamma_i = -1$ reduces exactly to this case. However, we will see that we need to exploit the freedom associated with $\gamma_i$ to optimize the local convergence properties of the difference map. We also see that only the compositions $\pi_i$ with $f_j$, for $i \neq j$, make sense when considered locally. This is because locally the constraint subspace $C_i$ may be approximated as an affine linear space, where $\pi_i \circ f_i = \pi_i$, and hence nothing is gained by introducing the maps $f_i$. Similarly, it is easily shown that nothing beyond the form (\ref{fMap}) is gained by considering arbitrary points on the {\it plane} determined by all three points ($\rho$ and its two projections, $\pi_1(\rho)$ and $\pi_2(\rho)$).

The optimum values for $\gamma_1$ and $\gamma_2$ are determined by considering one iteration of the difference map for points $\rho$ in the vicinity of a solution. Putting aside the case of histogram constraints, we assume that the constraint subspaces $C_1$ and $C_2$ are smooth, or nearly so. At a solution, $C_1$ and $C_2$ intersect and can be approximated as linear affine spaces. We will treat the case where the corresponding linear spaces $X_1 = C_1 - C_1$ and $X_2 = C_2 - C_2$ are orthogonal. Although this assumption may fail for certain pairs of projections, it prevails in the probabilistic sense of randomly oriented spaces satisfying $\dim{X_1}+\dim{X_2} < N$ in the limit $N \to \infty$. If $Y$ is the complement in $E^N$ of the span of $X_1$ and $X_2$, then a general point can be uniquely expressed as $\rho = x_1 + x_2 + y$, where $x_i \in X_i$ and $y \in Y$. The constraint spaces are now explicitly approximated as
\begin{equation}
\begin{array}{rl}
C_1 &= X_1 + a_2 + b_1 \\
C_2 &= a_1 + X_2 + b_2 \quad,
\end{array}
\end{equation}
where $a_i \in X_i$ and $b_i \in Y$. Orthogonality of the spaces $X_1$, $X_2$ and $Y$ implies the shortest element in $C_1 - C_2$ is $b_1 - b_2 \in Y$. Thus a solution (intersection) corresponds to $b_1 = b_2$, while $b_1 \neq b_2$ represents a ``trap", that is, a source of stagnation. The latter is illustrated by the action of the Gerchberg-Saxton map $G = \pi_1 \circ \pi_2$. First note the formulas for projections of a general point (these make use of the distance minimizing property of $\pi_1$ and $\pi_2$):
\begin{equation}
\begin{array}{rl}
\pi_1(x_1 + x_2 + y) &= x_1 + a_2 + b_1\\
\pi_2(x_1 + x_2 + y) &= a_1 + x_2 + b_2 \quad.
\end{array}
\end{equation}
One application of $G$ brings us to the fixed point
\begin{equation}
G(x_1 + x_2 + y) = a_1 + a_2 + b_1
\end{equation}
and stagnation occurs. Subsequent applications of the elementary projections simply hop between $a_1 + a_2 + b_1$ and $a_1 + a_2 + b_2$, the two points on $C_1$ and $C_2$ with minimum separation.

The behavior of the difference map is quite different, as we now show. A straightforward calculation gives the result
\begin{equation}\label{Daction}
D(x_1 + x_2 + y) = a_1 + a_2 + y + (1 - \beta \gamma_2)(x_1 - a_1) + (1 + \beta \gamma_1)(x_2 - a_2) + \beta(b_1 - b_2) \quad.
\end{equation}
First consider the case $b_1 = b_2 = b$, corresponding to a true intersection of the subspaces at the solution $\rho_{1\cap 2} = a_1 + a_2 + b$. From (\ref{Daction}) we see that subsequent iterates approach the fixed point $\rho^* = a_1 + a_2 + y$ provided $0 < \beta \gamma_2  <  2$ and $-2 < \beta \gamma_1 < 0$. This excludes $\gamma_1 = \gamma_2 = -1$, or identity maps, for the $f_i$. Optimal convergence (one iteration) is achieved by the parameter values $\gamma_2 = \beta^{-1}$ and $\gamma_1 = -\beta^{-1}$. This choice also makes $D$ maximally contractive; the rank of the Jacobian is then $\dim{Y}$. As explained in the previous Section, a fixed point of $D$ guarantees a true solution and it is given by $\rho_{1\cap 2} = (\pi_1 \circ f_2) ( \rho^* ) = a_1 + a_2 + b$.

Next we examine the situation for $b_1 \neq b_2$, when the two constraint subspaces locally form a trap. With the parameters $\gamma_i$ set to their optimal values, subsequent iterates have the form
\begin{equation}
D^n(x_1 + x_2 + y) = a_1 + a_2 + y + n \beta(b_1 - b_2) \quad.
\end{equation}
Rather than hopping between the two constraint subspaces we see that the iterates move away uniformly along the nearest separation axis, $b_1 - b_2$, with a rate determined by $\beta$. In a sense, the map has recognized that there is no solution at hand locally and seeks one elsewhere.

The local geometry of solutions and traps is rendered schematically in Figure 2. 
\begin{figure}[h]
\centerline{\scalebox{1.}{\includegraphics{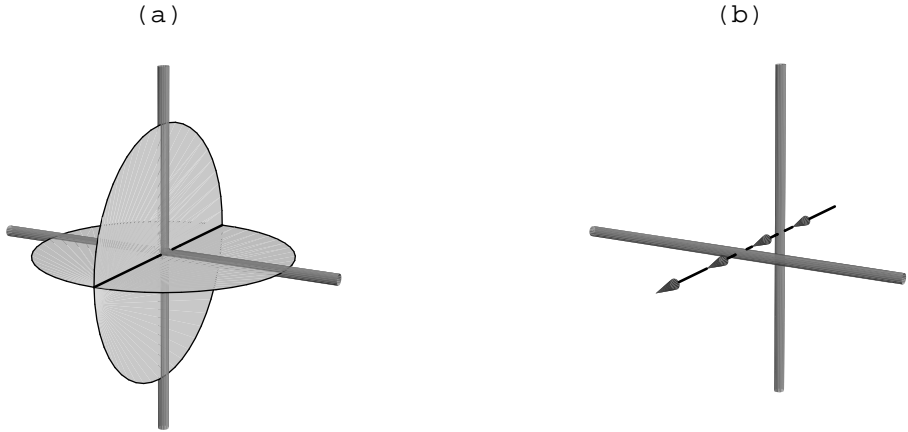}}}
\caption{Constraint subspaces (perpendicular rods) in the neighborhood of a solution (a) and a trap (b). The two circular disks in (a) represent points which project to the intersection of the constraint subspaces, $\rho_{1\cap 2}$, under action of $\pi_1 \circ f_2$ and $\pi_2 \circ f_1$ respectively; their intersection is the set of fixed points of the difference map. When the constraint subspaces do not intersect, as in (b), the action of the difference map is to move the iterates along the axis of minimum separation.
}
\end{figure}
It should be understood that this view, and the analysis in this Section, is valid only when there exists a minimum local separation of the two constraint subspaces, $\delta_{\rm min} = \| b_1 - b_2 \|$, and $\delta_{\rm min}$ is small on the scale where these subspaces deviate from linearity. We summarize the results of the last two Sections with the final form of the difference map:
\begin{equation}\label{diffMap2}
D: \rho \mapsto \rho + \beta \left( \pi_1\left[ (1 + \beta^{-1})\pi_2(\rho) - \beta^{-1} \rho \right] - \pi_2\left[ (1 - \beta^{-1})\pi_1(\rho) + \beta^{-1} \rho\right] \right) \quad.
\end{equation}
We note that all nonzero values of $\beta$ offer interesting choices, and that reversing the sign of $\beta$ has the effect of interchanging the two projections.

\section{Relationship to Fienup's input-output maps}

From (\ref{diffMap2}) we see that four elementary projections are computed for every iteration of the difference map. This number is reduced to two when $\beta = \pm 1$. In the case of support constraints, that is, $\pi_1 = \pi_{\rm supp}$ and $\pi_2 = \pi_{\rm mod}$, and the choice $\beta = 1$, one obtains the simple result:
\begin{equation}\label{hybrid}
D\big|_{\beta = 1} = F_{\rm hybrid}\; \colon\quad \rho_n \mapsto \rho^{\prime}_n = \left\{
\begin{array}{ll}
\pi_{\rm mod}(\rho)_n & \mbox{if $n \in S$}\\
\rho_n - \pi_{\rm mod}(\rho)_n &\mbox{if $n \notin S$} \quad.
\end{array}
\right.
\end{equation}
This is Fienup's hybrid input-output map\cite{Fienup1} for an object with support $S$ and Fienup's parameter $\beta_{\rm F} = 1$. In his discussion of input-output maps, Fienup also considered the maps
\begin{equation}
F_{\rm in-out} \colon\; \rho_n \mapsto \rho^{\prime}_n = \left\{
\begin{array}{ll}
\rho_n & \mbox{if $n \in S$}\\
\rho_n - \beta_{\rm F}\pi_{\rm mod}(\rho)_n &\mbox{if $n \notin S$} \quad.
\end{array}
\right.
\end{equation}
\begin{equation}
F_{\rm out-out} \colon\; \rho_n \mapsto \rho^{\prime}_n = \left\{
\begin{array}{ll}
\pi_{\rm mod}(\rho)_n & \mbox{if $n \in S$}\\
(1- \beta_{\rm F})\pi_{\rm mod}(\rho)_n &\mbox{if $n \notin S$} \quad.
\end{array}
\right.
\end{equation}
Only the hybrid map (\ref{hybrid}) can be obtained as a special case of the difference map. It is interesting that of the three input-output maps tested numerically by Fienup, the hybrid map outperformed the others and the optimum performance occurred for the parameter value $\beta_{\rm F} \approx 1$. In these studies\cite{Fienup1} a fixed number of iterations of an input-output map was followed by several iterations of the Gerchberg-Saxton map in order to arrive at the solution. In view of the discussion in Section 3, concerning the distinction between the solution and the fixed point of the map, the Gerchberg-Saxton iterations can be avoided for $\beta_{\rm F} = 1$. In this case the solution $\rho_{1\cap 2}$ is obtained from the fixed point $\rho^*$ by an elementary projection (in effect a single Gerchberg-Saxton iteration):
\begin{equation}
\rho_{1\cap 2} = (\pi_2 \circ f_1)( \rho^* ) = \pi_{\rm mod}( \rho^* ) \quad.
\end{equation}
No such statement can be made for other values of $\beta_{\rm F}$, for which there is no corresponding difference map.

The choice $\beta = -1$ in the difference map with support constraint is also relatively simple, but does not appear to have been studied previously. This has the effect of interchanging the roles of support and Fourier modulus projection. To express the resulting map most compactly we introduce a sign flipping operation on an object with support $S$:
\begin{equation}
R_S \colon\quad \rho_n \mapsto \rho^{\prime}_n = \left\{
\begin{array}{ll}
\rho_n & \mbox{if $n \in S$}\\
-\rho_n & \mbox{if $n \notin S$}\quad.
\end{array}
\right.
\end{equation}
The $\beta = -1$ counterpart of the hybrid input-output map is then given by
\begin{equation}
D\big|_{\beta = -1} \colon\quad \rho_n \mapsto \rho^{\prime}_n = \left\{
\begin{array}{ll}
(\pi_{\rm mod} \circ R_S)(\rho)_n & \mbox{if $n \in S$}\\
\rho_n + (\pi_{\rm mod} \circ R_S)(\rho)_n & \mbox{if $n \notin S$} \quad.
\end{array}
\right.
\end{equation}
In this case the solution is obtained from the fixed point most directly by
\begin{equation}
\rho_{1\cap 2} = (\pi_1 \circ f_2)( \rho^* ) = \pi_{\rm supp}( \rho^* ) \quad.
\end{equation}
Numerical experiments (Sec. 7) show that the difference map with $\beta = -1$ is often just as effective in finding solutions as the choice $\beta = 1$.

\section{Atomicity projection and Sayre's equation}

Phase retrieval in crystallography has traditionally been performed purely in the Fourier domain, whereas the {\it a priori} constraint of greatest relevance, atomicity, is most directly expressed in the object domain. A key development in the history of crystallographic phase retrieval was Sayre's equation\cite{Sayre}, which permitted a simple translation from the object to the Fourier domain for the case of identical, non-overlapping atoms. In the object domain Sayre's equation reads
\begin{equation}
\rho = g\ast ( \rho\times\rho ) \quad,
\end{equation}
where $\rho\times\rho$ denotes the nonlinear operation on the vector $\rho\in E^N$ which squares each component, and $\ast$ is the discrete convolution operator. The object $g$ is the point spread function corresponding to the finite atomic size $\sigma$. To simplify the analysis we consider the continuum limit, where $\sigma$ is large on the scale of the unit spacing of pixels. Using $r$ to denote pixel positions relative to some origin, then for a Gaussian atom in $d$ dimensions, the choice
\begin{equation}
g_{r} = \left( \frac{8}{\pi\sigma} \right)^{d/4}\exp{-\frac{2|r|^2}{\sigma}}
\end{equation}
is consistent with Sayre's equation when $\rho=\rho^{(1)}$, where
\begin{equation}
\rho^{(1)}_{r} = \left( \frac{2}{\pi\sigma} \right)^{d/4}\exp{-\frac{|r-r_0|^2}{\sigma}}\quad.
\end{equation}
We note that Sayre's equation ``selects" not just the size of atoms, $\sigma$, but also their normalization. The choice given above corresponds to
\begin{equation}\label{atomNorm}
{\| \rho \|}^2 = {{1}\over{M}} \sum_{r} {| \rho_{r} |}^2 =1 \quad,
\end{equation}
where $M$ is the number of atoms. Sayre's equation is easily translated into the Fourier domain by interchanging componentwise multiplication and convolution.

The earliest implementation of Sayre's equation in the Fourier domain was the ``tangent formula," a mapping of phases given by
\begin{equation}
T\colon\; \arg{\tilde{\rho}} \to \arg{\left[ \tilde{g}\times ( \tilde{\rho}\ast\tilde{\rho} ) \right]} \quad.
\end{equation}
Since the Fourier modulus of $\tilde{\rho}$ is fixed, $T$ is a mapping of the subspace $C_{\rm mod}$ onto itself. Solutions to the phase problem for identical atoms must be fixed points of $T$, but the converse need not be true since equality of the phases does not preclude inequality of the moduli:
\begin{equation}
|\tilde{\rho}| \neq \left| \tilde{g}\times ( \tilde{\rho}\ast\tilde{\rho} ) \right|\quad.
\end{equation}
In fact, a naive application of $T$ is unstable to the trivial ``uranium atom" solution, where all the phase angles are zero (or a translation thereof).

A more sophisticated treatment\cite{Debaerdemaeker} utilizes an objective function constructed from Sayre's equation:
\begin{equation}\label{SayreFunc}
V =\frac{1}{2}\, \| \rho - g\ast ( \rho\times\rho ) \|^2\quad.
\end{equation}
Given some small $\epsilon > 0$ which quantifies the error in the characterization of atomicity by Sayre's equation for the problem at hand, the region of $E^N$ satisfying $V<\epsilon$ can be identified with the subspace of atomicity constraints, $C_{\rm atom}$. With the gradient of $V$,
\begin{equation}\label{gradient}
\nabla V\colon\; \rho \to \rho - g\ast ( \rho\times\rho ) - 2(\bar{g}\ast \rho)\times\rho + 2\left( \bar{g}\ast g\ast(\rho\times\rho)\right)\times\rho\quad,
\end{equation}
where $\bar{g}_{r}=g_{-r}$, we can hope to find $C_{\rm atom}$ using the method of steepest descent. Within the usual framework of phase retrieval as practiced in crystallography, one would consider the map
\begin{equation}
S_{\rm mod}=\pi_{\rm mod}\circ  (1 - \alpha\nabla V)\quad.
\end{equation}
In the limit $\alpha\to 0^{+}$, the iterates of $S_{\rm mod}$ flow in the direction of decreasing $V$ on the subspace of modulus constraints, $C_{\rm mod}$. Apart for being impractical, the evolution for $\alpha\to 0^{+}$ is plagued by the problem of proliferating local minima (stagnation). Numerical experiments suggest that stagnation occurs over a range $0<\alpha <\alpha_{\rm c}$, where at the upper limit only the most attractive minima, including the desired global one, are active. Thus by carefully tuning $\alpha$ near $\alpha_{\rm c}$, the map $S_{\rm mod}$ can provide a reasonably practical phase retrieval scheme. A modification of the map for $\alpha=1$, which resembles the tangent formula when expressed in the Fourier domain,
\begin{equation}
\left. S_{\rm mod}\right|_{\alpha = 1}\; \colon \arg{\tilde{\rho}}\to \arg{\left[ \tilde{g}\times ( \tilde{\rho}\ast\tilde{\rho} ) + 2(\tilde{\bar{g}}\times \tilde{\rho})\ast\tilde{\rho} - 2\left( |\tilde{g}|^2\times (\tilde{\rho}\ast\tilde{\rho})\right)\ast\tilde{\rho}\right]}\quad,
\end{equation}
was used successfully in small molecule structure determination from x-ray data\cite{Debaerdemaeker}.

Sayre's equation can also be used in conjunction with the difference map, where, by treating Fourier modulus and atomicity constraints separately, stagnation is avoided. Rather than $S_{\rm mod}$, one considers
\begin{equation}
S_{\rm norm}=\pi_{\rm norm}\circ  (1 - \alpha\nabla V)\quad,
\end{equation}
where $\pi_{\rm norm}$ is the projection on the sphere with normalization (\ref{atomNorm}). The fixed points of $S_{\rm norm}$ are atomic objects, without regard to Fourier modulus, and thus can be identified with $C_{\rm atom}$. A good approximation to the projection operator is given by
\begin{equation}\label{SayreProj}
\pi_{\rm atom}\approx \pi_{\rm Sayre}= S_{\rm norm}^k
\end{equation}
for large $k$. The local distance minimizing properties of $\pi_{\rm Sayre}$, as required for  convergence of the difference map, can also be checked. Formally this follows from the smooth form of $\nabla V$, and the structure of its eigenvalues when linearized about a $V = 0$ fixed point. Since $V = 0$ is the global minimum, all such fixed points are attractive or equivalently, the linearization of $\nabla V$ has no negative eigenvalues. Zero eigenvalues will occur and correspond to continuous motion of the atoms, that is, motion along the subspace $C_{\rm atom}$. For sufficiently small $\alpha$, the linearization of $S_{\rm norm}$ is thus collapsing, on the eigenspaces with positive eigenvalues, and the identity on the local tangent space of $C_{\rm atom}$. By iterating, $S_{\rm norm}$ approaches a canonical (distance minimizing) projection operator.

A handle on the parameters $\alpha$ and $k$ is obtained from the stability analysis of a single atom having only a normalization degree of freedom: $\rho = \lambda \rho^{(1)}$.
Moreover, when the total number of atoms $M$ is large, we can neglect the action of $\pi_{\rm norm}$ and obtain $S_{\rm norm}(\lambda \rho^{(1)}) = f_{\alpha}(\lambda) \rho^{(1)}$, where
\begin{equation}\label{falpha}
f_{\alpha}\colon\;\lambda\to \lambda - \alpha\left[ \lambda - \lambda^2 + c(\lambda^3 - \lambda^2)\right]\quad,
\end{equation}
and $c=2(4/3)^{d/2}$. The map $f_{\alpha}$ has attractive fixed points $\lambda = 0, 1, \pm\infty$, where the latter is the ``uranium" instability. There is also a repulsive fixed point $\lambda^{\ast}=c^{-1}$ which separates the ``non-atom" (0) and ``atom" (1) fixed points. For small $\alpha$ the flow from $\lambda^{\ast}$ to either 0 or 1 is slow, and many iterations $k$ are needed. A reasonable criterion for selecting $k$ is that the slope of $f_{\alpha}^k$ at $\lambda^{\ast}$ is greater than one by a significant factor, say 2:
\begin{equation}\label{k_alpha}
\log{2}<\log{\left[f_{\alpha}^{\prime}(\lambda^{\ast})\right]^{k}}=k\log{\left[ 1+\alpha(1-c^{-1})\right]}\approx k\alpha(1-c^{-1})\quad.
\end{equation}
Figure 3 
\begin{figure}[h]
\centerline{\scalebox{1.}{\includegraphics{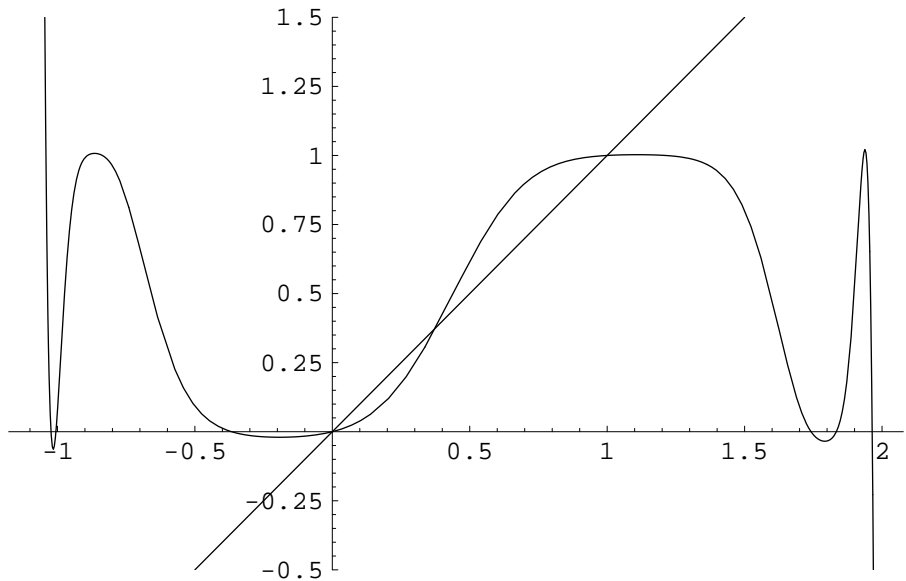}}}
\caption{Plot of the map $f_{\alpha}$, Eq. (\ref{falpha}), composed with itself three times for $\alpha=0.37$ and $c=8/3$. Points of intersection with the straight line give the three fixed points $\lambda=0,c^{-1},1$; the divergent behavior near $\lambda=-1$ and $\lambda=2$ corresponds to the ``uranium" instability.
}
\end{figure}
shows a plot of $f_{\alpha}^3$ with $\alpha = 0.37$, chosen so that (\ref{k_alpha}) is an equality for $d=2$ ($c=8/3$). The uranium instability sets in shortly beyond the upper boundary of the $\lambda=1$ fixed point's basin of attraction, given by the largest root of
\begin{equation}
0=\frac{f_{\alpha}(\lambda_{0})-1}{\lambda_{0}-1}=1-\alpha(c\lambda_{0}^2-\lambda_{0})\quad.
\end{equation}
If we view $\lambda_{0}$ as the largest starting value that does not encounter the uranium instability, then
\begin{equation}\label{alpha}
\alpha<\frac{1}{c\lambda_{0}^2 - \lambda_{0}}\quad.
\end{equation}
When the Sayre projection (\ref{SayreProj}) is applied to an object $\rho$ comprising a large number of atoms $M$, the starting normalizations $\lambda_{0}$ of the different atoms are drawn from a normal distribution. If the largest of these normalizations is in conflict with (\ref{alpha}) the projection will fail. The maximum of $M$ elements drawn from a normal distribution has the extreme value distribution\cite{extremeValue} with mean that grows as $\sqrt{\log{M}}$. Inequality (\ref{alpha}) thus implies that $\alpha$ be decreased roughly as the logarithm of the number of atoms.

With appropriately chosen $k$ and $\alpha$, the difference map implementation of Sayre's equation does not suffer from stagnation, as is common with tangent formula based schemes. Comparisons with the more direct approach to atomicity projection, as described in Section 2, are given in the next Section.  

\section{Numerical experiments}

Apart from the case of support constraints in one dimension\cite{nonunique1d}, the role of dimensionality in phase retrieval seems to be minor. This is the conclusion of experiments performed in one, two and three dimensions using the difference map with histogram and atomicity constraints. The length of computations (difference map iterations) required to locate $M$ atoms in the crystal unit cell, or $M$ stars in an image of the sky, or $M$ spikes in a time series, are empirically very comparable, given equal numbers of Fourier moduli in the appropriate dimension. Other properties of the object, for example compactness, appear to be far more important in determining the computational complexity of phase retrieval. Given this ambivalence, most of the experiments below are two dimensional in conformity with the dimensionality of the print medium.

A uniform set of normalization and initialization conventions were used in all the experiments. All objects $\rho$ (points in $E^N$) were normalized so that $\| \rho \| = 1$.
When the Fourier modulus and histogram is normalized with the same convention, normalization is not required during the course of the iterations.
\begin{figure}[p]
\centerline{\scalebox{1.}{\includegraphics{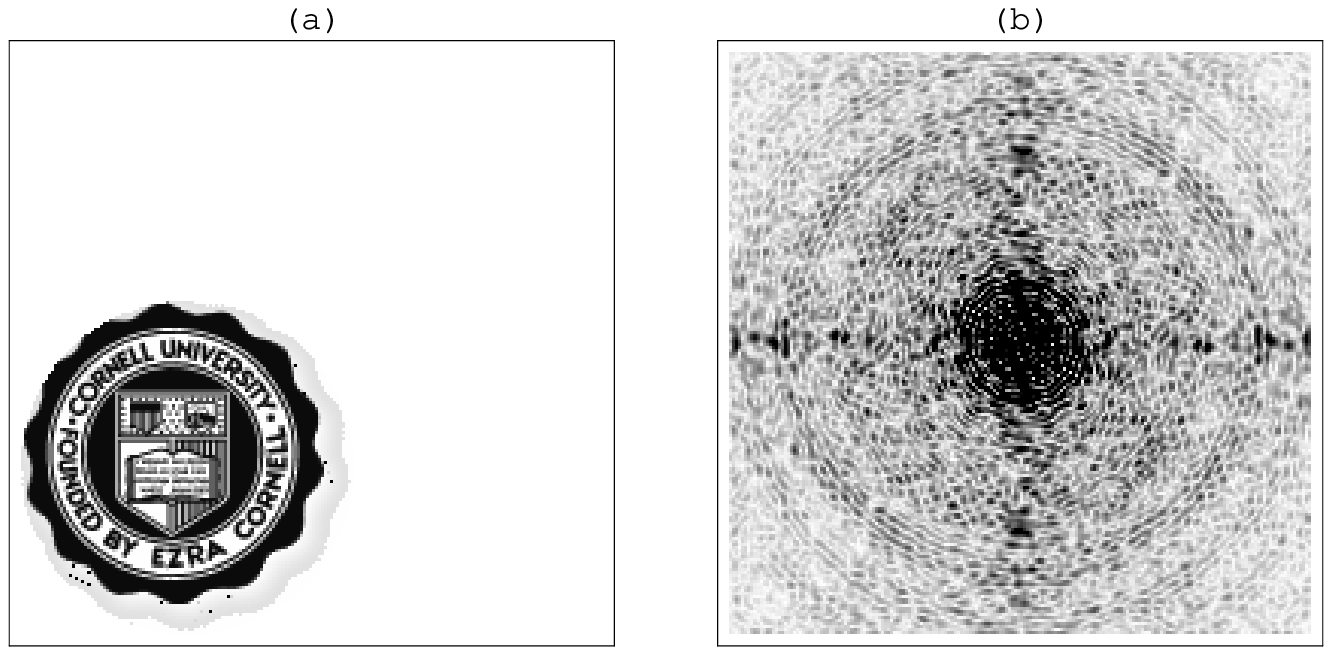}}}
\caption{Object (a) and corresponding Fourier modulus (b) used in subsequent phase retrieval experiments.
}
\centerline{\scalebox{1.}{\includegraphics{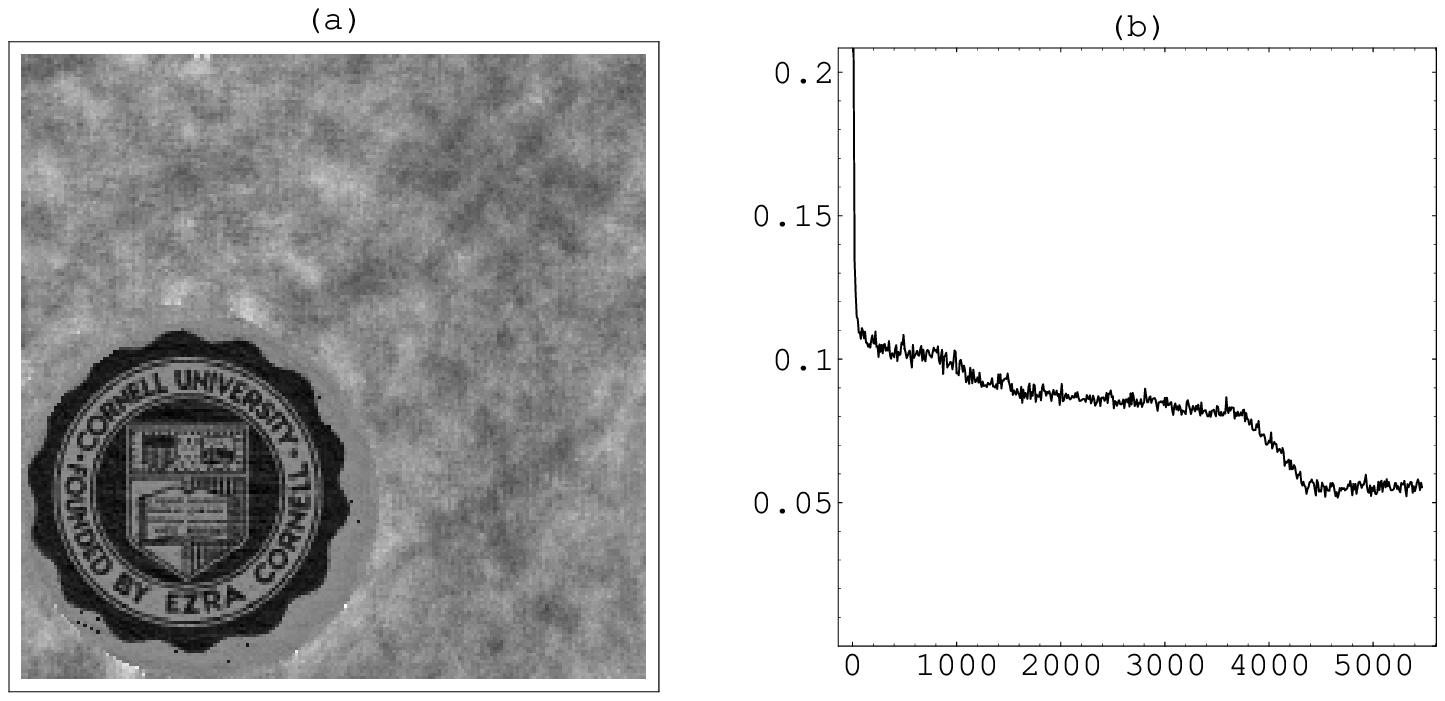}}}
\caption{Phase retrieval using the difference map with support and positivity constraints. (a) Final iterate (fixed point of the map); (b) behavior of the error $e_i$ with iteration $i$. The random gray contrast in (a) is removed, and the original object (Fig. 4(a)) is almost perfectly recovered, by a single application of Fourier modulus projection.
}
\end{figure}
With support and positivity constraints, normalization was applied after every projection. The components of the initial object $\rho(0)$ were generated by a pseudo-random number generator and then normalized. All the experiments used Fourier modulus projection as the projection $\pi_2$.

Figure 4 shows the $192\times 192$ pixel image and corresponding Fourier modulus used in the first experiments. The difference map is first demonstrated using support and positivity constraints ($\pi_1 = \pi_{\rm supp}\circ\pi_{\rm pos}$) and the value $\beta = 1$, for which it reproduces Fienup's hybrid input-output map. A circle, measuring 112 pixels in diameter, was used as support. After only about 10 iterations from the random start an object comprising a circular ring of dark contrast was formed. Subsequently, as seen in the series of plateaus in the error plot shown in Figure 5(b),
the various characteristics of the object were refined. The 14-fold modulation of the outer ring appears after about 1000 iterations; the text is revealed after  another 3000 iterations, etc. After 5000 iterations the magnitude of the difference between successive iterates (as given quantitatively in the error plot) is so small that the object has arrived at the essentially static fixed point  shown in Figure 5(a). The best estimate of the solution, $\rho_{1\cap 2}$, is obtained by a single application of $\pi_{\rm mod}$ to the last iterate and is almost indistinguishable from the true object (Fig. 4(a)). We note that the random gray contrast in the fixed point image (Fig. 5(a)) is a manifestation of the fact that fixed points of the difference map live in a large space, the particular point chosen being subject to the randomness in $\rho(0)$. Very similar results (not shown) are obtained with $\beta = -1$, whereas performance degrades both for small and large $\beta$ of either sign.

Histogram constraints are much more effective in retrieving the phases of the object in Figure 4(a). With $\pi_1 = \pi_{\rm hist}$ and $\beta = 1$, the object is perfectly recovered in only about 100 iterations, as shown in Figure 6.
\begin{figure}[t]
\centerline{\scalebox{1.}{\includegraphics{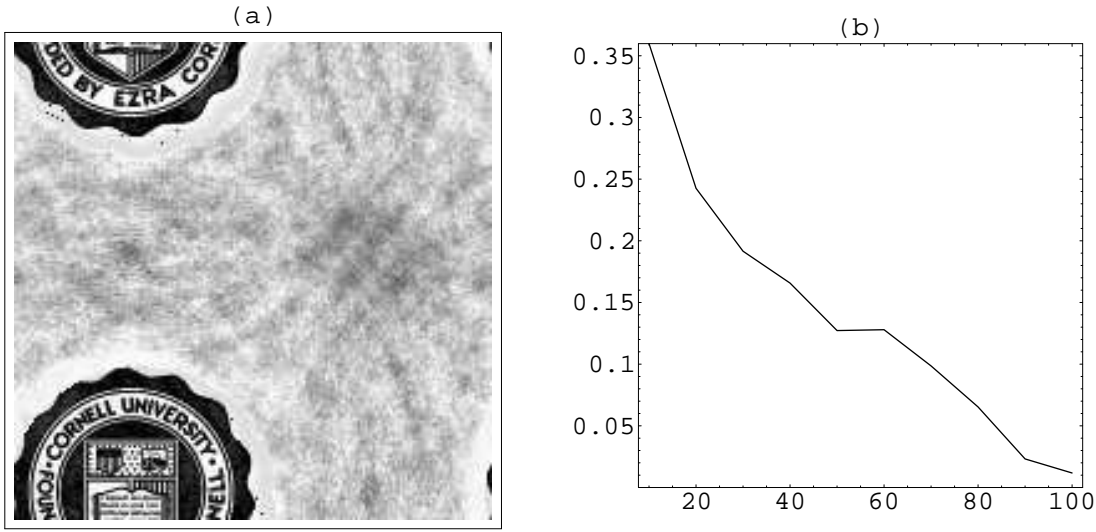}}}
\caption{Same as Figure 5 but with histogram constraints. Far fewer iterations are required to recover the same object.
}
\end{figure}
Again, the mottled background in the fixed point object, Figure 6(a), is completely eliminated by a single Fourier modulus projection. Results with $\beta = -1$ are again very similar, although there are qualitative differences in the appearance of the iterates. The larger initial error using histogram projection also indicates that the corresponding constraint subspace is significantly smaller than the subspace of support constraints used in the previous experiment. We conclude that for this particular object the histogram appears to contain more information than knowledge of the object support.

The ``histogram" (sorted list of $192\times 192$ object values) used to implement histogram projection for the previous example is shown plotted in Figure 7(a). 
\begin{figure}[t]
\centerline{\scalebox{1.}{\includegraphics{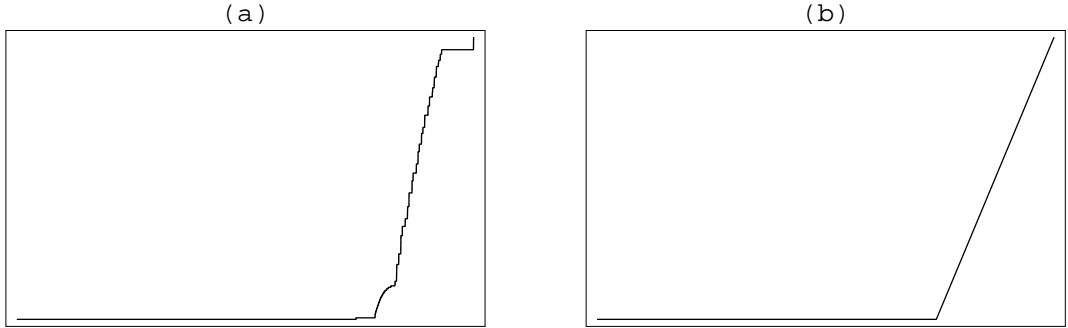}}}
\caption{Plots of sorted pixel values, or ``histograms." (a) Cornell seal, Fig. 4(a);(b) random disk, Fig. 8(a). 
}
\end{figure}
Because of the quirks in the histogram, peculiar to this particular object, the histogram data cannot be viewed as a valid form of {\it a priori} knowledge. An example that comes closer to a valid application of histogram constraints is the object shown in Figure 8(a). 
\begin{figure}[t]
\centerline{\scalebox{1.}{\includegraphics{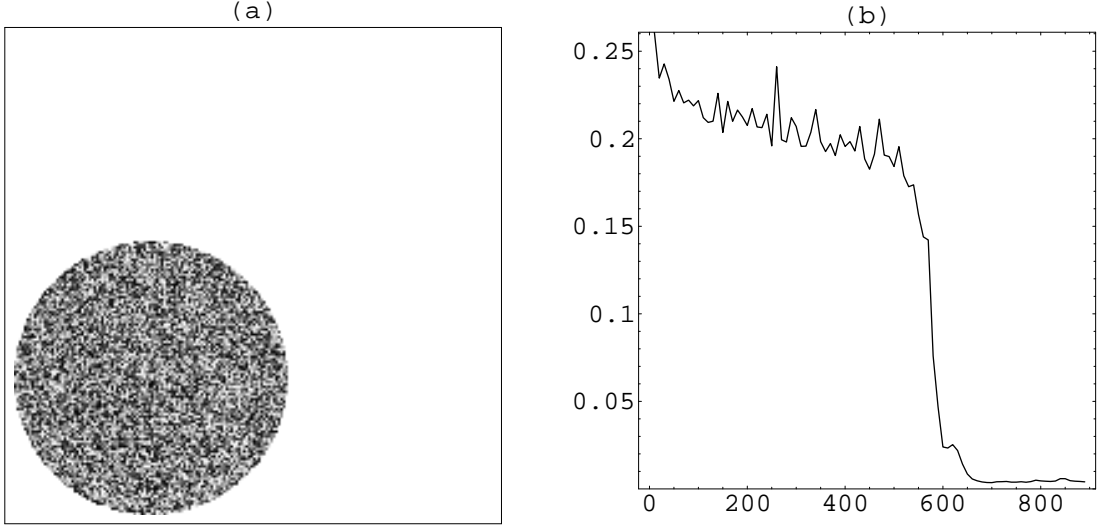}}}
\caption{Phase retrieval applied to a random disk (a) using histogram constraints; (b) plot of the error. The recovered object (not shown) not only had the correct circular outline but reproduced in detail the pixel values within the disk. 
}
\end{figure}
This object has the same circular support as the previous object. However, its pixel values were assigned random numbers from a uniform distribution; the resulting histogram is shown in Figure 7(b). Although the progress toward the solution was slower initially, the detailed pattern of random pixel values of this object was perfectly recovered after about 600 iterations (Fig. 8(b)).

Histogram constraints apply ideally in situations where the object is composed of individual elements, e.g. atoms, whose values are known in detail even if the location of each element is not. If we refer to objects with this property generically as ``atomic", then the histogram of atomic objects is specified uniquely by the number and distribution of their constituent ``atoms". For simplicity, all the experiments below were performed with $M$ identical atoms, each atom being an approximation to a Gaussian on a small number of pixels. In the two dimensional experiments the support of each atom was a $3\times 3$ array of pixels; the continuously variable center of the atom with respect to the central pixel then gives rise to a continuous distribution of values for the histogram. Details regarding the construction of approximate Gaussians for a general choice of atomic support are given in the Appendix.

Objects having a specified number of atoms $M$ were generated by applying atomicity projection to a random image,
\begin{equation}
\rho_{\rm atom} = \pi_{\rm atom}(\rho_{\rm rand}) \quad.
\end{equation}
To simulate the effects of clustering, additional power was optionally applied to the low spatial frequency components of $\rho_{\rm rand}$. Specifically, the Fourier modulus of $\rho_{\rm rand}$ was specified using the formula
\begin{equation}\label{avePower}
{| \tilde{\rho}_{\rm rand} |^2}_q = {{q_0^2}\over{|q|^2 + q_0^2}} \quad,
\end{equation}
where $|q|$ is the magnitude of the Fourier wavevector and $q_0 = 2\pi/\xi$ is the characteristic wavevector corresponding to a clustering length scale $\xi$. The limit $\xi \to 0$ reproduces a white power spectrum and uncorrelated atomic positions in $\rho_{\rm atom}$ (apart from avoided overlaps). For finite $\xi$ the atoms in $\rho_{\rm atom}$ are clustered in groups with linear dimension of order $\xi$. The phases of $\tilde{\rho}_{\rm rand}$  were drawn at random from a uniform distribution, giving the maximally random atomic object possible. We note that the Fourier modulus after atomicity projection, $| \tilde{\rho}_{\rm atom} |$, is no longer a smooth function of $q$ although its behavior at small $|q|$ resembles Eq. (\ref{avePower}) when locally averaged over wavevector. At large $|q|$ the rapid decay of $| \tilde{\rho}_{\rm atom} |$ with $|q|$ reflects the approximate Gaussian Fourier transform of each individual atom.

Figure 9(a) 
\begin{figure}[p]
\centerline{\scalebox{1.}{\includegraphics{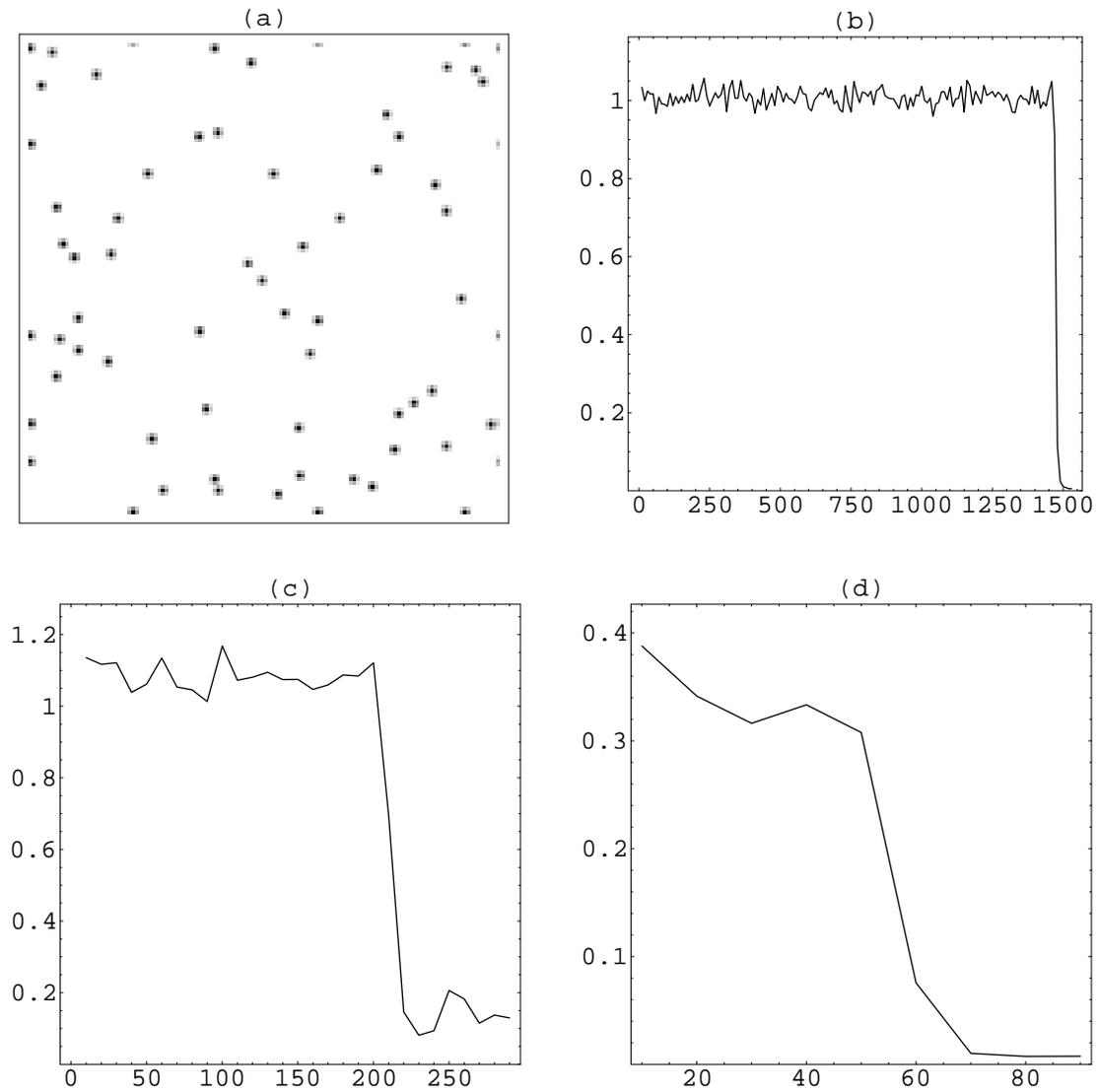}}}
\caption{(a) Object comprising 60 equal ``atoms"; error plots for three different implementations of {\it a priori} information: (b) histogram constraint, (c) $3\times 3$ pixel atoms, (d) atomicity implemented by Sayre's equation. 
}
\end{figure}
shows an example of 60 uncorrelated atoms ($\xi \to 0$) on $128\times 128$ pixels. The adjoining error plots document the successful phase retrieval for this object by difference map iterations with three different choices for $\pi_1$: (b) histogram projection and $\beta = 1$, (c) atomicity projection with $3\times 3$ pixel Gaussian atoms and $\beta = 0.5$, (d) atomicity projection using Sayre's equation ($\pi_1 = \pi_{\rm Sayre}$, Eq. (\ref{SayreProj}), $k=3$, $\alpha=0.37$) and $\beta = 1$. The error plots suggest that the phases are retrieved by a more nearly random process than is typically the case for compact objects, such as the ones considered earlier. A large number of experiments support the view that with uncorrelated atoms the search for the solution is a deterministic but nevertheless {\it random} walk through a relatively homogeneous ``landscape". The step size in this walk, as measured by the angle $\theta_{1\cap2}$, is very nearly constant and only shrinks significantly when, apparently by accident, the walk arrives at the attractive basin of the solution. That there is no manifest progress toward the solution until the very end of the process is confirmed by snapshots of the object along the way. In contrast, the reconstruction using histogram information of the compact objects in the previous experiments proceeded through definite stages of refinement. Although the Sayre projection, Fig. 9(d), found the solution relatively quickly, experiments with larger numbers of atoms displayed the same random walk behavior exhibited by the other projections. We also note that because the Sayre projection is much more demanding computationally than atomicity projection in the object domain ($\pi_{\rm atom}$), the performance of the two forms of atomicity projection is not as different as the error plots imply. Object domain atomicity projection has the added advantage that multiple kinds of atoms are easily accommodated.

The random walk behavior of the difference map for uncorrelated atoms is also observed with other values of $\beta$. A quantitative comparison that is free of object idiosyncrasies is achieved by averaging results over the ensemble of atoms defined above. The results of a study with uncorrelated ($\xi \to 0$) atoms on $64\times 64$ pixel images is shown in Figure 10. 
\begin{figure}[t]
\centerline{\scalebox{1.}{\includegraphics{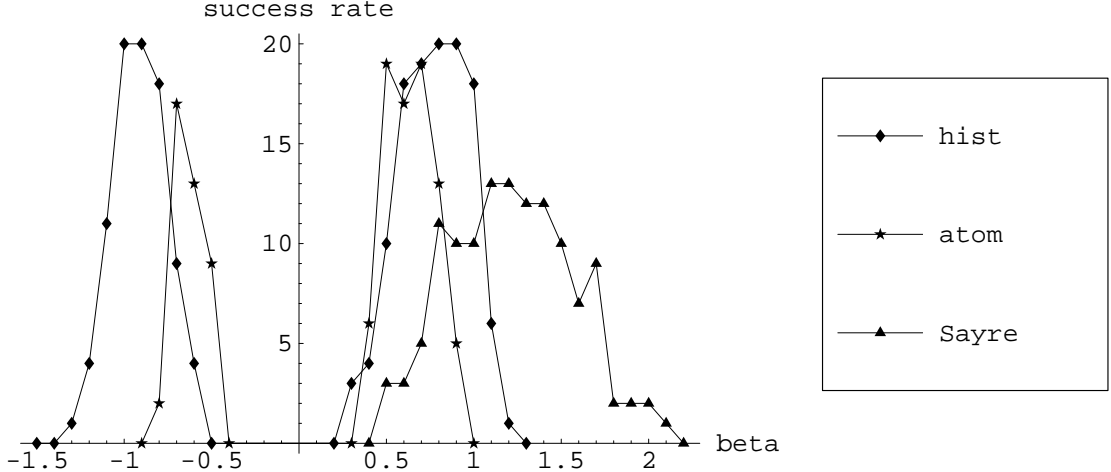}}}
\caption{Comparison of phase retrieval successes after 100 iterations of the difference map for histogram projection, $3\times 3$ pixel atomicity projection, and  atomicity projection via Sayre's equation. The horizontal axis is the difference map parameter $\beta$. See text for more details.
}
\end{figure}
For each value of $\beta$ shown, and three choices for $\pi_1$ ($\pi_{\rm hist}$, $\pi_{\rm atom}$ and $\pi_{\rm Sayre}$), 20 experiments were performed, each with a different realization of atoms and 100 iterations of the difference map. In the studies of $\pi_{\rm hist}$ and $\pi_{\rm atom}$ systems of 30 atoms were used and the success rate has peaks at both positive and negative $\beta$. Since for this number of atoms $\pi_{\rm Sayre}$ gave a 100\% success rate over a large range of $\beta$, the number of atoms was increased to 40. The success rate for this projection has a broader peak at larger $\beta$ and virtually zero success at negative $\beta$ due to the ``uranium" instability. We note that $\pi_{\rm Sayre}$ has the additional parameters $k$ and $\alpha$ that can be tuned to improve performance. This study used $k=3$ and $\alpha=0.37$.

The number of difference map iterations required for phase retrieval grows rapidly with the number of atoms and is the subject of a future study. Here we make the observation that any degree of clustering in the atomic positions has a profound effect, apparently reducing the number of iterations by orders of magnitude. Figure 11(a) 
\begin{figure}[t]
\centerline{\scalebox{1.}{\includegraphics{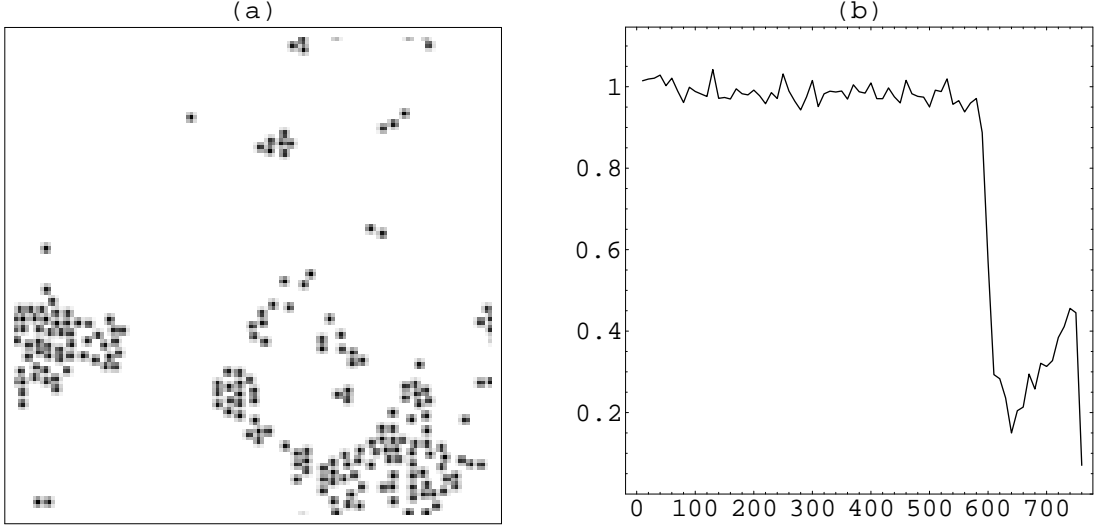}}}
\caption{Phase retrieval becomes easier when the atoms are clustered, as in this example (a) comprising 200 atoms; error plot for the case of $3\times 3$ pixel atomicity projection.
}
\end{figure}
shows an example with 200 atoms generated with a clustering length scale $\xi=50$ on $128\times 128$ pixels. Without clustering, a system of this many atoms could not be solved using the difference map in a reasonable time. The effect of the clustering is to greatly reduce the search space, in a manner not unlike what was observed for the compact objects (Figs. 4(a) and 8(a)) studied earlier, where enhanced power in the low spatial frequencies helps locate the ``support" of the clusters. Snapshots of the iterates show a relatively slow variation of low frequency features accompanied by rapid fluctuations in the atom positions. The example shown was solved in about 600 iterations using $\pi_1=\pi_{\rm atom}$ and $\beta = 0.5$. The fluctuations in the error (Fig. 11(b)), after the solution has been found, correspond to random translations of the entire object.

For completeness we include two experiments with atomic objects in dimensions other than two. As already mentioned, the performance of the difference map was found to be indistinguishable from studies in two dimensions having the same total number of pixels and atoms, all other attributes (clustering, atomic size, etc.) being the same. Figure 12 
\begin{figure}[t]
\centerline{\scalebox{1.}{\includegraphics{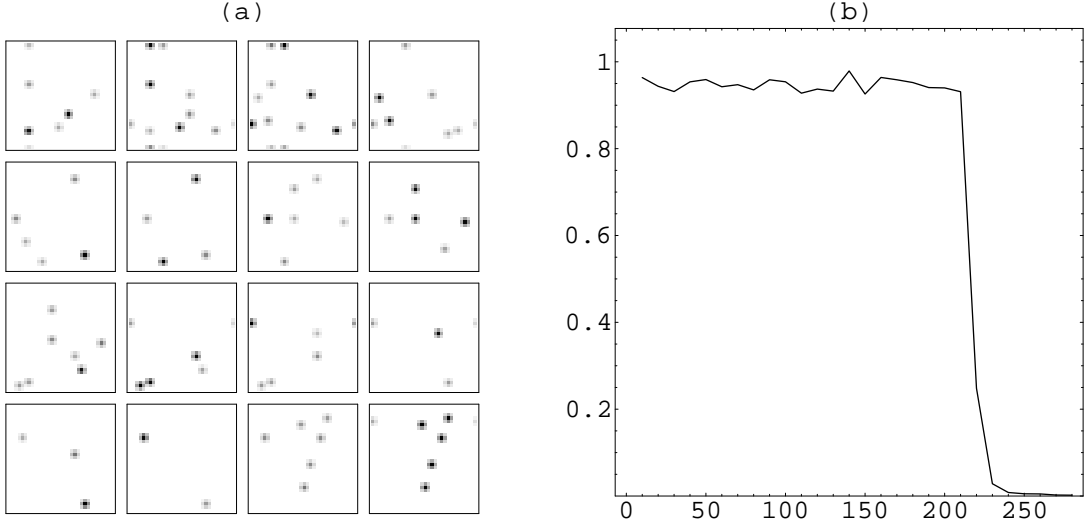}}}
\caption{Phase retrieval for an object comprising 60 atoms in three dimensions. (a) The first 16 layers, arranged lexicographically, of the $32\times 32\times 32$ voxel array; (b) error plot for the difference map with histogram projection.
}
\end{figure}
shows an experiment with 60 uncorrelated $3\times 3\times 3$ voxel atoms on a $32\times 32\times 32$ grid; 16 layers of the object are arranged in Fig. 12(a). The difference map using histogram projection and $\beta =1$ found the solution in about 250 iterations. An experiment in one dimension is documented in Figure 13. 
\begin{figure}[p]
\centerline{\scalebox{1.}{\includegraphics{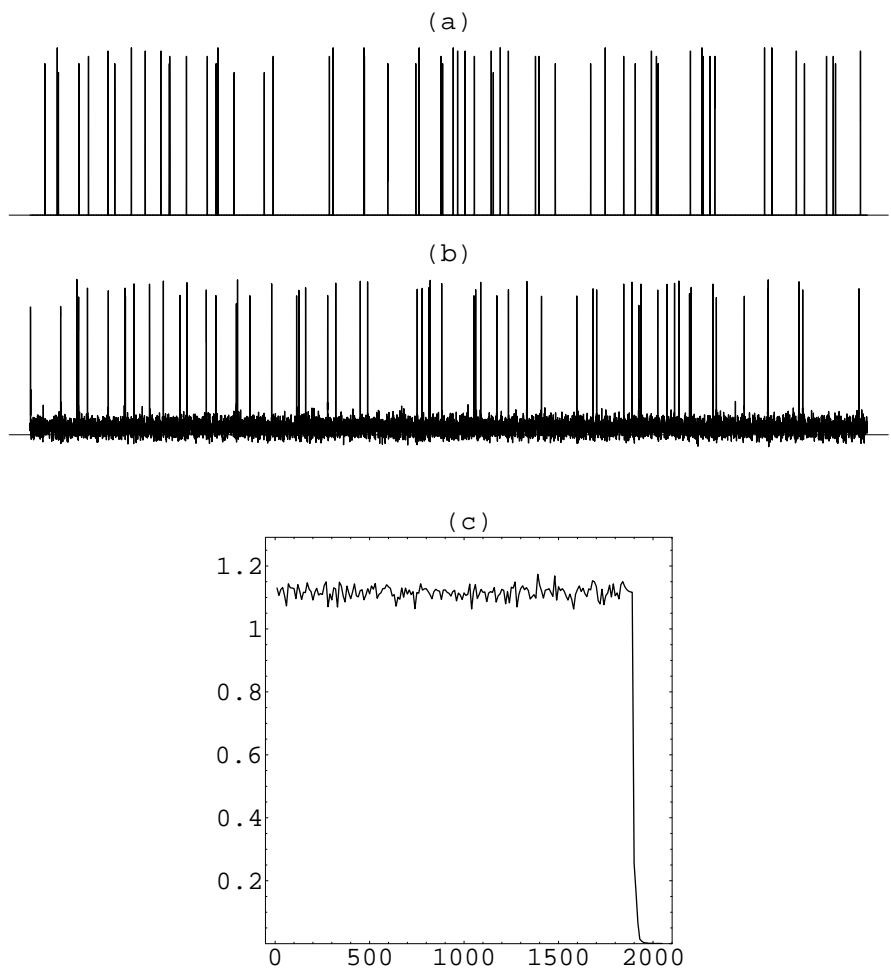}}}
\caption{Phase retrieval for an object comprising 60 ``atoms" in one dimension. (a) Plot of object values; (b) fixed-point of difference map with histogram projection; (c) error plot. A single application of Fourier modulus projection to (b) perfectly reproduces (a), up to inversion and translation (not shown).
}
\end{figure}
The plot of the object values, Fig. 13(a), shows 60 uncorrelated ``atoms", each with support on $3$ adjacent grid points; the size of the object is $2^{14}=16384$, the same as the two dimensional example in Figure 9(a). Figure 13(b) shows the fixed point of the difference map, again with histogram projection and $\beta =1$. One application of Fourier modulus projection to the object in Fig. 13(b) perfectly reproduces a translated and inverted copy of the original object, Fig. 13(a). We see that the number of iterations required to find the solution (Fig. 13(c)), about 2000, is very similar to the result obtained in the two dimensional experiment, Fig. 9(b).

\section{Conclusions}
We are still far from being able to claim that an algorithm for phase retrieval is at hand. Like other iterative schemes currently in use, the difference map considered here falls short of a true algorithm in that there are no useful bounds on the number of iterations required to find the solution (a fixed point). The analysis given in Section 3 only guarantees that every fixed point corresponds to a solution (by Eq. (\ref{projFixedPt})). Fixed points of dynamical systems are just the simplest form of an {\it attractor}; deterministic systems are also known to have ``strange attractors", complex spaces upon which iterates of the map can become trapped. From numerical experiments there is some evidence that the existence of strange attractors, in competition with the simple fixed points, is the main failure mechanism (stagnation) of the difference map when $\beta$ is too small. For small $\beta$ the iterates exhibit chaotic fluctuations while preserving certain gross features in the object. Strange attractors can be tolerated when their number is small: one need only restart the iterations, hopefully within the  basin of attraction of a simple fixed point. This is not an option when strange attractors proliferate, apparently for small $\beta$, and possibly also in situations where the phase retrieval problem is close to being underdetermined.

Competing strange attractors appear to be absent when $\beta$ is not small and the phase retrieval problem is well posed (overdetermined). In this regime the numerical evidence suggests that the solution will always be found given enough iterations. By construction (Sec. 4), the difference map is highly contractive and we are led to believe that the fraction of ``phase space" actually explored by the map is extremely small. 

Some insights about the nature of the space being explored can be gained from the study of phase retrieval experiments where the data has been fabricated to remove all fixed points (solutions) while preserving the general characteristics of the map. One such fabrication, for example, is to form two different objects having the same histogram and attempt phase retrieval using their averaged Fourier moduli. It is highly unlikely that phases exist which can be combined with this fabricated modulus data and yield an object having the required histogram. On the other hand, because the statistical properties of the fabricated data is very similar to that of either of the two genuine objects, the qualitative behavior of the difference map should be unchanged. The only real difference between the genuine and fabricated data is that the difference map iterates for the former will eventually find a fixed point, while for the latter the iterates will continue exploring indefinitely a space that we expect to share all the essential characteristics of the space being explored in the case of the genuine objects. This space is itself a strange attractor and the complexity of iterative phase retrieval is perhaps best quantified by some measure of its size.

A ``phase portrait" of a strange attractor, obtained with fabricated data of the kind described above, is shown for three values of $\beta$ in Figure 14. 
\begin{figure}[t]
\centerline{\scalebox{1.}{\includegraphics{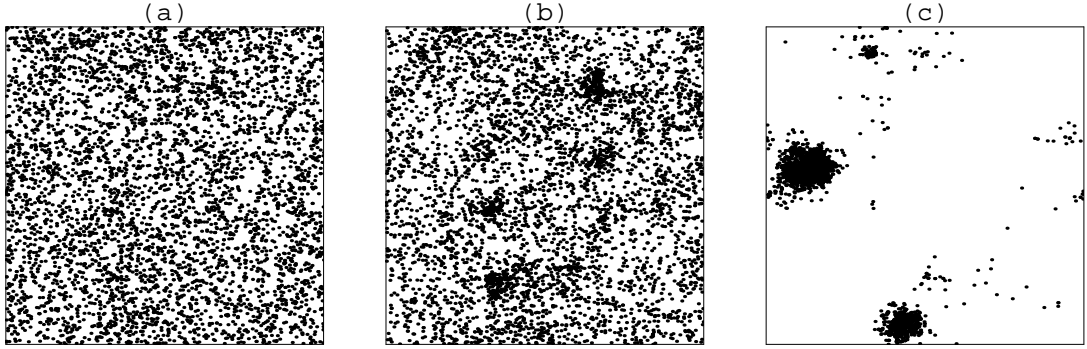}}}
\caption{Evolution of the difference map attractor with decreasing $\beta$: (a) 0.7, (b) 0.6, (c) 0.5.
}
\end{figure}
The two genuine objects were one dimensional with identical two-valued histograms: distinct sequences of 16 zeros and 16 ones. Iterates were obtained for the difference map with histogram projection. Two dimensional images of the attractor were generated by first projecting each iterate into the four-torus corresponding to the first four (unrestricted) phase angles of the Fourier transform. Selecting only those iterates whose second pair of phase angles were both close to zero, the first pair of angles were plotted as a point in the plane. The images are thus projections of codimension-2 sections of the attractor. Figure 14 shows the evolution of the attractor from a relatively uniform probability distribution for $\beta=0.7$ to a highly nonuniform distribution for $\beta=0.5$. The collapse of the search space with decreasing $\beta$ is consistent with a corresponding improvement in phase retrieval performance (for non-fabricated data). The absence of inversion symmetry in the image for $\beta=0.5$ indicates that the threshold to the regime of multiple attractors has been crossed.

The study of iterative phase retrieval in the context of discrete dynamical systems may provide new insights for improving solution strategies and even a handle on a comprehensive theory of the complexity of algorithms. A meaningful measure of the size of the search space, for example, might be given in terms of the dimensionality of the associated attractor. It would be interesting to study not just the evolution of complexity (e.g. attractor dimensionality) with parameters defining the iterative map ($\beta$), but also with respect to quantifiable attributes of the object (clustering). 

The application of phase retrieval which currently poses the greatest challenge is macromolecular crystallography. Recently the Shake and Bake (SnB) method\cite{SnB2} has succeeded in retrieving phases for small proteins without the benefit of the additional MIR or MAD data normally required for protein structure determination. SnB\cite{SnB1} is iterative and uses object domain atomicity projection and objective function minimization as its two elementary operations. On general grounds, however, one can argue that the SnB method may not have realized its full potential. First, the two operations are implemented in the alternating Gerchberg-Saxton fashion which in most other applications leads to stagnation. Second, the objective function used by SnB is based on an uncorrelated collection of identical point scatterers, not unlike the function $V$ (Eq. (\ref{SayreFunc})) derived from Sayre's equation. An objective function (or associated projection) that makes use of clustering, or the division of space into solvent and non-solvent regions, would add considerable {\it a priori} information and thereby reduce the size of the space being searched. An example of the improved performance with tighter object domain constraints was seen in the experiment above which compared histogram and atomicity projection (Figs. 9(b) and 9(c), respectively).

\section{Appendix}

\subsection{Finitely sampled Gaussians}

An ``atom" in $d$-dimensions with center $r_0$ is modelled as the Gaussian
\begin{equation}
\Psi(r,r_0) = \left(\frac{2}{\pi \sigma}\right)^{d/4}\exp{-\frac{| r - r_0 |^2}{\sigma}} \quad,
\end{equation}
with normalization
\begin{equation}
\int | \Psi |^2 d^d r = 1 \quad.
\end{equation}
We are interested in approximating $\Psi$ on a small set of pixels $S = \{ s_1, s_2, \dots \}$ on a $d$-dimensional cubic grid with unit spacing. For a given choice of ``atomic support" $S$, there is an optimum Gaussian width $\sigma$ that we aim to determine.

The support $S$ is defined in a translationally invariant way relative to the atom center $r_0$, or rather, the grid point $p_0$ closest to $r_0$. For randomly placed atoms we expect the fractional translation $t = r_0 - p_0$ to be distributed uniformly in the cube $T = (-\frac{1}{2},\frac{1}{2})^d$.

Let $\{ \tilde{\Psi}_s(t) \mid s \in S \}$ be the approximation on $S$ of the Gaussian $\Psi$ with fractional center $t \in T$. The normalization
\begin{equation}\label{norm}
\sum_{s \in S} | \tilde{\Psi}_s(t) |^2 = 1
\end{equation}
applies for each $t \in T$. Also, for each $t \in T$, the values $\tilde{\Psi}_s(t)$ are determined by minimizing the Euclidean distance between $\tilde{\Psi}$ and the restriction of $\Psi$ to $S$:
\begin{equation}
\delta(t) = \sum_{s \in S} | \tilde{\Psi}_s(t) - \Psi(s,p_0 + t) |^2 \quad.
\end{equation}
Minimizing $\delta(t)$ subject to Eq. (\ref{norm}) gives
\begin{equation}
\tilde{\Psi}_s(t) = \frac{\Psi(s,p_0 + t)}{\sqrt{\sum_{s \in S} | \Psi(s,p_0 + t) |^2}} \quad,
\end{equation}
and the squared Euclidean distance
\begin{equation}
\delta(t) = \left( \sqrt{\sum_{s \in S} | \Psi(s,p_0 + t) |^2} -1\right)^2 \quad.
\end{equation}
Finally, the optimal Gaussian width $\sigma$ is determined by minimizing the average of $\delta(t)$ over fractional translations:
\begin{equation}
\delta_{\rm ave} = \int_T \delta(t) d^d t \quad.
\end{equation}
The minimization of $\delta_{\rm ave}$ with respect to $\sigma$ can be carried out numerically on a suitably fine grid of fractional translations $t$. Results for various choices of support in $d = 1$, $2$ and $3$ are given in Table 1. A natural choice for $S$ is the set of pixels within a distance $R=1, \sqrt{2},\sqrt{3},\dots$ of the origin. We expect the width $\sigma$ to increase and the distance $\delta_{\rm ave}$ to decrease as $R$ increases.

\begin{table}[t]
\centering
{\renewcommand{\arraystretch}{1}
\begin{tabular}
{|c||c|c|c|c|}
\hline
$d$&$R$& pixels &$\sigma$&$\delta_{\rm ave}$\\ \hline
$1$&$1$&$3$&$1.156$&$0.000025$\\

     &$2$&$5$&$1.800$&$4\times 10^{-8}$\\

     &$3$&$7$&$2.445$&$6\times 10^{-11}$\\ \hline

$2$&$1$&$5$&$0.814$&$0.0030$\\

     &$\sqrt{2}$&$9$&$1.115$&$0.000060$\\

     &$2$&$13$&$1.238$&$0.000021$\\ \hline

$3$&$1$&$7$&$0.694$&$0.011$\\

     &$\sqrt{2}$&$19$&$0.952$&$0.00070$\\

     &$\sqrt{3}$&$27$&$1.091$&$0.00010$\\ 
\hline
\end{tabular}
}
\caption{Widths $\sigma$ and errors $\delta_{\rm ave}$ of finitely sampled Gaussians in dimension $d$ for supports within a distance $R$ of the origin.}
\end{table}

The order of choosing the atomic support and corresponding Gaussian width $\sigma$ is usually reversed in actual applications. For example, in crystallography the Gaussian decay of intensity with diffraction wavevector  $q$ (due to atomic size, disorder, etc.) fixes the value of $\sigma$. For $M$ identical Gaussian atoms at positions $r_m$ we have 
\begin{equation}
\tilde{\rho}_q = \left(\frac{\sigma}{2\pi}\right)^{d/4} A(q) \exp{-\frac{\sigma}{4}|q|^2} \quad,
\end{equation}
where
\begin{equation}
A(q) = \sum_{m=1}^M \exp{i q\cdot r_m}
\end{equation}
has essentially a white spectrum for random $r_m$. A good estimate of $\sigma$ is therefore given directly in terms of the expectation value of $|q|^2$ with respect to the intensity data $|\tilde{\rho}_q|^2$:
\begin{equation}
\sigma^{-1} = \frac{1}{D}\frac{\sum_q |q|^2 | \tilde{\rho}_q |^2}{\sum_q | \tilde{\rho}_q |^2} \quad.
\end{equation}
After determining $\sigma$ one would then consult a table such as Table 1 to discover the atomic support $S$ for which this $\sigma$ is optimal, or close to optimal.

\subsection{Atomicity projection}

We assume that finitely sampled Gaussians $\{ \tilde{\Psi}_s(t) \mid s \in S \}$ have been determined for the problem at hand. These are precomputed on the chosen support $S$ with a suitably fine grid of fractional translations $t$. An object comprising $M$ identical atoms will have support $S+\{ p_1, p_2,\dots, p_M \}$, where the $p_m$ are atomic support centers on the integer grid. Normally one imposes a non-overlapping condition such as
\begin{equation}\label{overlap}
p_m - p_n \notin S - S \quad (m \neq n).
\end{equation}

Atomicity projection of an arbitrary (real valued) object $\rho$ is accomplished in three steps:

\noindent {\bf Obtain support centers.} To identify positions in $\rho$ having a large overlap with a Gaussian atom, $g$, one forms the convolution $\rho^{\prime}=g\ast\rho$. The pixel values of $\rho^{\prime}$ are then sorted and their locations are appended, beginning with the largest, to the list of atomic support centers $p_m$. Each potential new center must satisfy a separation condition such as (\ref{overlap}).

\noindent {\bf Obtain fractional positions.} Within each atomic support $S+p_0$ the actual atomic position may have a fractional translation $t$ as described above. Since the finitely sampled Gaussian $\{ \tilde{\Psi}_s(t) \mid s \in S+p_0 \}$ is an approximation to the true Gaussian $\Psi(r, p_0+t)$, its centroid (on $S+p_0$) will be close to the true Gaussian center $p_0+t$. To find $t$ one therefore obtains the centroid of $\rho^{\prime}$ restricted to $S+p_0$. If $\rho^{\prime}$ is far from being atomic it may happen that the $t$ determined via the centroid is outside the cube $T$ of allowed fractional translations. It is then necessary to obtain the translation $t^{\prime}$ on the boundary of $T$ that is closest to $t$.

\noindent {\bf Object synthesis} The final step is to add up the individual atomic objects. If atomic supports are allowed to overlap, the result must be normalized.

Although this projection algorithm is probably not distance minimizing for a general object $\rho$, it is nearly so in the situation which matters most, that is, when $\rho$ is already nearly atomic. We recall that the distance minimizing property of projections was required only in the analysis of local convergence (Sec. 4).

\section{Acknowledgment}
I thank John Spence for inviting me to the Workshop on New Approaches to the Phase Problem for Non-Periodic Objects (Lawrence Berkeley National Laboratory, May 2001), where, perhaps for the first time, phase retrieval experts from several different disciplines were assembled under one roof. This work was supported by the National Science Foundation under grant ITR-0081775.

\end{document}